\documentclass[12pt]{amsart}

\makeatletter

\newcommand{\LyX}{L\kern-.1667em\lower.25em\hbox{Y}\kern-.125emX\spacefactor1000}
\newcommand{\noun}[1]{\textsc{#1}}

\theoremstyle{plain}    
\newtheorem{thm}{Theorem}[section]
\numberwithin{equation}{section} 
\numberwithin{figure}{section} 
\theoremstyle{plain}    
\newtheorem{cor}[thm]{Corollary} 
\theoremstyle{plain}    
\newtheorem{prop}[thm]{Proposition} 
\theoremstyle{definition}
\newtheorem{defn}[thm]{Definition}
\theoremstyle{remark}
\newtheorem{rem}[thm]{Remark}
\theoremstyle{remark}    
\newtheorem{notation}[thm]{Notation} 

\newcommand{\Tr}{{\operatorname{Tr}}}
\newcommand{\sa}{{\operatorname{sa}}}

\makeatother

\begin{document}

\title{A Microstates Approach to Relative Free Entropy}

\author{Dimitri Shlyakhtenko}

\date{\today}

\thanks{Research supported in part by a National Science Foundation postdoctoral fellowship.}

\email{shlyakht@member.ams.org}

\address{Department of Mathematics, UCLA, Los Angeles, CA 90095}

\maketitle

\section{Introduction.}

\newcommand{\grb}[6]{\Gamma _{R}(#1 |#2 [#3 ];#4 ,#5 ,#6 )}
 
\newcommand{\Xn}{X_{1},\dots ,X_{n}}

\newcommand{\gr}[5]{\Gamma _{R}(#1 ;#2 ,#3 ,#4 )}
 
\newcommand{\Ym}{Y_{1},\dots ,Y_{m}}

\newcommand{\ym}{y_{1},\dots ,y_{m}}
 
\newcommand{\xn}{x_{1},\dots ,x_{n}}

\newcommand{\Zr}{Z_{1},\dots ,Z_{r}}

Let \( \Xn \in (M,\tau ) \) be a family of non-commutative random variables
in a tracial \( W^{*} \)-probability space, and let \( B\subset M \) be a
unital subalgebra. Voiculescu has introduced in \cite{dvv:entropy5} a free
entropy quantity
\[
\chi ^{*}(\Xn :B).\]
 His approach involved non-commutative Hilbert transform and is algebraic in
nature. In the case that \( B=\mathbb C \), this quantity is denoted \( \chi ^{*}(\Xn ) \),
and its properties are very similar to those of the free entropy \( \chi (\Xn ) \)
introduced by Voiculescu in \cite{dvv:entropy2} using microstates; in fact,
it may very well be that the two quantities coinside. 

Using the microstates approach to free entropy, we introduce in this paper a
quantity \( \chi (\Xn |B) \) which has several properties in common \( \chi ^{*}(\Xn :B) \). 

The infinitesimal change of variables formula for \( \chi  \) involves conjugate
variables, introduced by Voiculescu in order to define \( \chi ^{*} \). Both
\( \chi  \) and \( \chi ^{*} \) have the same behavior under compression by
matrix units in the case that \( B=M_{n}\otimes D \), where \( M_{n} \) is
the algebra of \( n\times n \) matrices. This behavior is a useful technical
tool; for example, it was used to prove certain maximization results for matrices
of non-commutative random variables in \cite{nss:entropy-micro}.

If the \( \Xn  \) are free from \( B \), we have
\[
\chi (\Xn |B)=\chi (\Xn ),\]
provided that \( B \) can be embedded into the ultrapower of the hyperfinite
II\( _{1} \) factor; a similar fact holds for \( \chi ^{*} \). 

If \( \Xn  \) are free from the algebra generated by \( \Ym  \) and \( B \),
then
\[
\chi (\Xn ,\Ym |B)=\chi (\Xn |B)+\chi (\Ym |B);\]
a similar fact holds for \( \chi ^{*} \).

We prove a maximization result for \( \chi  \) (which is essentially identical
to the one for \( \chi ^{*} \)), namely \( \chi (\Xn |B) \) attains its maximum
among all \( \Xn  \) with \( \sum _{i}\tau (X_{i}^{2})=n \) if and only if
\( \Xn  \) is a free semicircular family, free from \( B \) (we need as an
assumption that \( B \) can be embedded into the ultrapower of the hyperfinite
II\( _{1} \) factor). Lastly, the infinitesimal change of variables formula
for \( \chi (\cdots |B) \) involves conjugate variables used to define \( \chi ^{*}(\cdots :B) \)
(see \cite{dvv:entropy5}).

It is interesting to note that \( \chi (\cdot |B) \) has an interpretation
as a relative entropy, which suggests a similar interpretation for \( \chi ^{*}(\cdot :B) \).
Indeed, we show that if \( \Ym  \) are generators of \( B \), then \( \chi (\Xn |B)=\chi (\Xn |\Ym ) \),
where the latter entropy has properties of a relative entropy of \( \Xn  \)
and \( \Ym  \). We caution the reader that we use a definition of \( \chi (\Xn |\Ym ) \)
which may be different from the one used by Voiculescu in \cite{dvv:entropy2},
although the two quantities are related.

\section{Relative Free Entropy \protect\( \chi (\cdot |B)\protect \).}

Let \( (M,\tau ) \) be a tracial non-commutative probability space, and consider
self-adjoint non-commutative random variables \( X_{1},\dots ,X_{n} \), \( Y_{1},\dots ,Y_{m}\in M \).
We denote by \( M_{k} \) the algebra of \( k\times k \) matrices, and by \( M_{k}^{\sa } \)
the set of self-adjoint \( k\times k \) matrices. Recall that the set
\[
\Gamma _{R}(\Xn ;k,l,\epsilon )\subset (M_{k}^{\sa })^{n}\]
was defined by Voiculescu in \cite{dvv:entropy2} as the set of those \( (\xn )\in (M^{\sa }_{k})^{n} \),
for which \( \Vert x_{i}\Vert \leq R \) and for any \( p\leq l \), and all
\( i_{1},\dots ,i_{p} \)
\[
|\tau _{n}(x_{i_{1}},\dots ,x_{i_{p}})-\tau (X_{i_{1}}\dots X_{i_{p}})|<\epsilon .\]
Here \( \tau _{n} \) stands for the normalized trace on the matrices (so that
\( \tau _{n}(1)=1 \)).

\begin{defn}
For \( y_{1},\dots ,y_{n}\in M_{k} \), define
\begin{eqnarray*}
\grb{\Xn }{\ym }{\Ym }{k}{l}{\epsilon } &  & \\
=\{\xn \in (M_{k}^{\sa })^{n}:(\xn ,\ym ) &  & \\
\in \gr{\Xn ,\Ym }{k}{l}{\epsilon }{}\} &  & 
\end{eqnarray*}
 
\end{defn}
\begin{rem}
\label{rmq: non-empty}Note that for this set to be nonempty, we must have that
\[
(\ym )\in \gr{(\Ym :\Xn }{l}{k}{\epsilon }{}.\]
This set would be empty if \( W^{*}(\Ym ) \) were not  embeddable into the
ultrapower of the hyperfinite II\( _{1} \) factor.
\end{rem}
\begin{defn}
Let \( \lambda  \) denote Lebesgue measure on \( (M_{k}^{\sa })^{n} \) corresponding
to its Hilbert space structure coming from the non-normalized trace. Define
\emph{}successively\emph{
\begin{eqnarray}
\chi _{R}(\Xn |\Ym ,l,\epsilon )= &  & \label{ean:defofentropy} \\
\lim _{k\to \infty }\frac{1}{k^{2}}\sup _{\begin{array}{c}
(\ym )\in \\
\Gamma _{R}(\Ym ;k,l,\epsilon )
\end{array}} &  & \\
\log \lambda \grb{\Xn }{\ym }{\Ym }{k}{l}{\epsilon }+\frac{n}{2}\log k & 
\end{eqnarray}
}
\end{defn}

\[
\chi _{R}(\Xn |\Ym )=\inf _{l,\epsilon }\chi _{R}(\Xn |\Ym ;l,\epsilon )\]

\[
\chi (\Xn |\Ym )=\sup _{R}\chi _{R}(\Xn |\Ym ).\]
 The last quantity is called the relative free entropy of the \( n \)-tuple
\( (\Xn ) \) with respect to the \( m \)-tuple \( (\Ym ) \).

If \( \omega  \) is a free ultrafilter on \( \mathbb N \), then one can also
define \( \chi ^{\omega }(\Xn |\Ym ) \) exactly as in (\ref{ean:defofentropy}),
but replacing \( \limsup  \) by \( \lim _{k\to \omega } \). 

\begin{rem}
\label{rmq:relative}It is not clear whether our definition of \( \chi (\Xn |\Ym ) \)
coninsides with that of Voiculescu (see \cite{dvv:entropy2}). His definition
corresponds to defining
\begin{eqnarray*}
\chi _{R}(\Xn |\Ym ;l,k,\epsilon )= &  & \\
\limsup _{k}\log \frac{\lambda \Gamma _{R}(\Xn ,\Ym ;k,l,\epsilon )}{\lambda \Gamma _{R}(\Ym ;k,l,\epsilon )}+\frac{n}{2}\log k. &  & 
\end{eqnarray*}
The connection to our definition can be made as follows: let
\[
f(\ym )=\lambda \grb{\Xn }{\Ym }{\ym }{k}{l}{\epsilon }\]
Then Voiculescu's definition corresponds to taking as \( \chi _{R}(\Xn |\Ym ;l,k,\epsilon ) \)
the average of \( f \) over \( \Gamma _{R}(\Ym ;l,k,\epsilon ) \). It follows
that the quantity obtained in our definition is bigger than that of Voiculescu.
We mention that it is possible to define, in the spirit of \cite{dvv:entropy2}
the relative entropy as
\[
\chi '(\Xn |\Ym )=\chi (\Xn ,\Ym )-\chi (\Ym :\Xn ).\]
Such a definition corresponds to defining \( \chi _{R}(\Xn |\Ym ;l,k,\epsilon ) \)
as the average of \( f \) over \( \Gamma _{R}(\Ym :\Xn ) \). We don't know
whether \( \chi ' \) coincides with Voiculescu's or our definition of \( \chi , \)
and whether Voiculescu's and our definitions are the same or different. Note,
however, that we always have:
\end{rem}
\begin{prop}
Let \( \Xn  \), \( \Ym  \) be non-commutative random variables. Then we have
\begin{eqnarray*}
\chi (\Xn ,\Ym )-\chi (\Ym )\leq  &  & \\
\chi (\Xn ,\Ym )-\chi (\Ym :\Xn )\leq  &  & \\
\chi (\Xn |\Ym ). &  & 
\end{eqnarray*}
In particular, if \( \chi (\Xn |\Ym )=-\infty  \), then \( \chi (\Xn ,\Ym )=-\infty  \).
\end{prop}
\begin{proof}
In the case that \( \chi (\Ym ) \) is finite, the inequalities follow from
the discussion in Remark \ref{rmq:relative}. If
\[
\chi (\Ym )=-\infty ,\]
then
\[
\chi (\Xn ,\Ym )=-\infty .\]
If \( \chi (\Ym )\neq -\infty  \), then
\[
\chi (\Xn |\Ym )=-\infty \]
implies that
\[
\chi (\Xn ,\Ym )\leq -\infty +\chi (\Ym )=-\infty .\]

\end{proof}
\begin{notation}
We shall write
\[
\Gamma _{R}(\Xn |\ym ;l,k,\epsilon )\]
for
\[
\grb{\Xn }{\ym }{\Ym }{l}{k}{\epsilon }\]
when the \( \Ym  \) are understood.
\end{notation}
\begin{prop}
If \( p<m \), then \( \chi (\Xn |\Ym )\leq \chi (\Xn |Y_{1},\dots ,Y_{p}) \),
and a similar inequality holds for \( \chi _{R} \) and \( \chi ^{\omega } \).
\end{prop}
\begin{proof}
We have the inclusion
\[
\Gamma _{R}(\Xn |\ym [\Ym ];l,k,\epsilon )\subset \grb{\Xn }{y_{1},\dots ,y_{p}}{Y_{1},\dots ,Y_{p}}{l}{k}{\epsilon }\]
for all \( (\ym )\in M_{k} \). It follows that
\begin{eqnarray*}
\sup _{(\ym )}\lambda \Gamma _{R}(\Xn |\ym [\Ym ];l,k,\epsilon )\leq  &  & \\
\sup _{(y_{1},\dots ,y_{p})}\lambda \grb{\Xn }{y_{1},\dots ,y_{p}}{Y_{1},\dots ,Y_{p}}{l}{k}{\epsilon } &  & 
\end{eqnarray*}
which implies the desired inequality.
\end{proof}
\begin{prop}
\( \chi (\Xn |\Ym )\leq \chi (\Xn :\Ym ) \).
\end{prop}
\begin{proof}
We clearly have
\[
\Gamma _{R}(\Xn |\ym [\Ym ];l,k,\epsilon )\subset \pi \Gamma _{r}(\Xn ,\Ym ;l,k,\epsilon ),\]
where \( \pi  \) denotes the projection from \( (M_{k}^{\sa })^{n}\times (M^{\sa })^{m} \)
onto \( (M_{k}^{\sa })^{n} \).
\end{proof}
\begin{cor}
\label{corr: max}Let \( c^{2}=\frac{1}{n}\sum \tau (X_{i}^{2}) \). Then
\[
\chi (\Xn |\Ym )\leq \frac{n}{2}\log 2\pi ec^{2}.\]
The same estimate holds for \( \chi _{R} \) and \( \chi ^{\omega }. \) 
\end{cor}
\begin{proof}
We have \( \chi (\Xn ,\Ym )\leq \chi (\Xn )\leq \frac{n}{2}\log 2\pi ec^{2} \),
the last inequality by \cite{dvv:entropy2}. 
\end{proof}
\begin{defn}
Let \( B\subset M \) be a unital subalgebra of \( M \). We define the free
entropy of \( (\Xn ) \) relative to \( B \) to be
\[
\chi (\Xn |B)=\inf _{m}\inf _{\Ym \in B}\chi (\Xn |\Ym ).\]

\end{defn}
If \( \omega  \) is a free ultrafilter on the natural numbers, then we define
\( \chi ^{\omega }(\Xn |B) \) in the obvious way. 

\begin{rem}
For \( \chi (\Xn |B) \) to be finite, we must have that for all \( \Ym \in B \),
\[
\chi (\Xn |\Ym )\neq -\infty .\]
By Remark~\ref{rmq: non-empty}, this subsumes that \( B \) is embeddable into
the ultrapower of the hyperfinite II\( _{1} \) factor.
\end{rem}
\begin{prop}
If \( D\subset B \) is a unital subalgebra, then
\[
\chi (\Xn |B)\leq \chi (\Xn |D).\]
In particular, we have \( \chi (\Xn |B)\leq \chi (\Xn |\mathbb {C})=\chi (\Xn ) \).
The same conclusion holds for \( \chi  \) replaced with \( \chi ^{\omega } \).
\end{prop}
\begin{proof}
The first inequality is because in computing \( \chi (\Xn |D) \) we take the
infimum over a smaller set. The equality between \( \chi (\Xn |\mathbb {C}) \)
and \( \chi (\Xn ) \) is left to the reader.
\end{proof}
\begin{prop}
Let \( X_{i}^{j} \), \( i=1,\dots ,n \), \( j=1,2,\dots  \) be non-commutative
random variables. Assume that \( \Xn  \) are such that as \( j\to \infty  \),
the joint distribution of \( (X_{1}^{j},\dots ,X_{n}^{j}) \) and \( B \) converges
to the joint distribution of \( B \) and \( \Xn  \). Let \( Y_{1}^{j},\dots ,Y_{m}^{j} \)
be such that they converge in distribution to \( \Ym  \). Assume that there
exists a finite constant \( R \), so that \( \sup _{i,j}\Vert X_{i}^{j}\Vert ,\Vert Y_{i}^{j}\Vert ,\Vert X_{i}\Vert ,\Vert Y_{i}\Vert <R \).
Then
\[
\chi (\Xn |B)\geq \limsup _{j}\chi (X_{1}^{j},\dots ,X_{n}^{j}|B)\]
and
\[
\chi (\Xn |\Ym )\geq \limsup _{j}\chi (X_{1}^{j},\dots ,X_{n}^{j}|Y_{1}^{j},\dots ,Y_{m}^{j}).\]
The same conclusion holds for \( \chi ^{\omega } \) instead of \( \chi  \).
\end{prop}
\begin{proof}
Clearly, only the second inequality needs to be proved. It follows from the
following inclusion, true for sufficiently large \( j \):
\begin{eqnarray*}
\Gamma _{R}(X_{1}^{j},\dots ,X_{n}^{j}|y_{1},\dots ,y_{m}[Y_{1}^{j},\dots ,Y_{n}^{j}];k,l,\epsilon /2)\subset  &  & \\
\Gamma _{R}(\Xn |\ym [\Ym ;k,l,\epsilon ), &  & 
\end{eqnarray*}
since it implies that
\begin{eqnarray*}
\sup _{\ym }\log \lambda \Gamma _{R}(X_{1}^{j},\dots ,X_{n}^{j}|y_{1},\dots ,y_{m}[Y_{1}^{j},\dots ,Y_{n}^{j}];k,l,\epsilon /2)\leq  &  & \\
\sup _{\ym }\log \lambda \Gamma _{R}(\Xn |\ym [\Ym ;k,l,\epsilon ), &  & 
\end{eqnarray*}
which in turn implies the desired inequality.
\end{proof}
\begin{prop}
\label{prop: generation}Assume that \( \Zr \in W^{*}(\Ym ) \). Then
\[
\chi (\Xn |\Ym )\leq \chi (\Xn |\Zr ),\]
and similarly for \( \chi  \) replaced with \( \chi ^{\omega } \).
\end{prop}
\begin{proof}
Fix \( \epsilon >0 \), \( l>0 \), \( R>0 \). Then there exists non-commutative
polynomials \( p_{1},\dots ,p_{r} \) in \( m \) variables, such that
\[
\left( \Xn ,p_{1}(\Ym ),\dots ,p_{r}(\Ym )\right) \]
approximate
\[
(\Xn ,\Zr )\]
strongly to any desired accuracy. It follows, that for a suitable choice of
such polynomials, there exist \( \epsilon >\epsilon '>0 \), \( l'>l>0 \) and
\( R'>R>0 \) such that whenever
\[
(\ym )\in \Gamma _{R}(\Ym ,\Xn ;k,l',\epsilon ')\]
we have the inclusion
\begin{eqnarray*}
\Gamma _{R'}(\Xn |\ym [\Ym ];k,l',\epsilon ')\subset  &  & \\
\Gamma _{R}(\Xn |p_{1}(\ym ),\dots ,p_{r}(\ym )[\Ym ];k,l,\epsilon ). &  & 
\end{eqnarray*}
It follows that
\begin{eqnarray*}
\sup _{\ym }\log \lambda \Gamma _{R'}(\Xn |\ym [\Ym ];k,l',\epsilon ')\leq  &  & \\
\sup _{z_{1},\dots ,z_{r}}\log \lambda \Gamma _{R}(\Xn |z_{1},\dots ,z_{r}[\Zr ];k,l,\epsilon ). &  & 
\end{eqnarray*}
But this implies \( \chi (\Xn |\Ym )\leq \chi (\Xn |\Zr ) \).
\end{proof}
\begin{thm}
\( \chi (\Xn |\Ym )=\chi (\Xn |W^{*}(\Ym )) \).
\end{thm}
\begin{proof}
By definition,
\[
\chi =\chi (\Xn |W^{*}(\Ym ))=\inf _{r}\inf _{\Zr \in W^{*}(\Ym )}\chi (\Xn |\Zr ).\]
In particular, if \( r=m \) and \( (\Zr )=(\Ym ) \), we have that
\[
\chi \leq \chi (\Xn |\Ym ).\]
But by Proposition \ref{prop: generation}, we also have that, for \( Z_{1},\dots ,Z_{r}\in W^{*}(\Ym ) \),
\[
\chi (\Xn |\Ym )\leq \chi (\Xn |\Zr )\leq \chi ,\]
so that \( \chi =\chi (\Xn |\Ym ) \).
\end{proof}
\begin{prop}
If \( \Ym \in B \), then \( \chi (\Xn |B)\leq \chi (\Xn :\Ym ) \).
\end{prop}
\begin{proof}
We have \( \chi (\Xn |B)\leq \chi (\Xn |\Ym )\leq \chi (\Xn :\Ym ) \).
\end{proof}
\begin{thm}
Let \( B=L^{\infty }[0,1]\subset M \) be a diffuse commutative von Neumann
subalgebra. For each \( r=1,\dots ,n \) let \( \mu _{r} \) be a measure on
\( [0,1]^{2} \) determined by 
\[
\iint f(x)g(y)d\mu _{r}(x,y)=\tau (X_{r}fX_{r}g).\]
 Assume that for at least one \( r\in \{1,\dots ,n\} \), Lebesgue measure on
\( [0,1]^{2} \) is singular with respect to \( \mu _{r} \). Then \( \chi (\Xn |B)=-\infty  \). 
\end{thm}
\begin{proof}
Let \( Y \) be a self-adjoint generator for \( B \). Then we have 
\[
\chi (\Xn |B)=\chi (\Xn |Y)\leq \chi (\Xn :Y)=-\infty \]
because of \cite[Corollary 7.7]{dvv:entropy3}. 
\end{proof}
Note that the analogous theorem holds for \( \chi ^{*}(\cdots :B) \), see \cite{shlyakht:cpentropy}. 

\begin{prop}
Let \( p<n \). Then
\[
\chi (\Xn |B)\leq \chi (X_{1},\dots ,X_{p}|B)+\chi (X_{m+1},\dots ,X_{p}|B)\]
and a similar inequality holds for \( \chi ^{\omega } \).
\end{prop}
\begin{proof}
We clearly have
\begin{eqnarray*}
\Gamma _{R}(\Xn |\ym ;k,l,\epsilon )\subset \Gamma _{R}(X_{1},\dots ,X_{p}|\ym ;k,l,\epsilon )\times  &  & \\
\Gamma _{R}(X_{p+1},\dots ,X_{n}|\ym ;k,l,\epsilon ) &  & 
\end{eqnarray*}
so that
\begin{eqnarray*}
\sup _{\ym }\log \lambda \Gamma _{R}(\Xn |\ym ;k,l,\epsilon )\leq  &  & \\
\sup _{\ym }\log \lambda \left( \Gamma _{R}(X_{1},\dots ,X_{p}|\ym ;k,l,\epsilon )+\right.  &  & \\
\left. \log \lambda \Gamma _{R}(X_{p+1},\dots ,X_{n}|\ym ;k,l,\epsilon \right) \leq  &  & \\
\sup _{\ym }\log \lambda \left( \Gamma _{R}(X_{1},\dots ,X_{p}|\ym ;k,l,\epsilon )\right) + &  & \\
\sup _{\ym }\log \lambda \left( \Gamma _{R}(X_{p+1},\dots ,X_{nb}|\ym ;k,l,\epsilon \right) . &  & 
\end{eqnarray*}
This implies the proposition.
\end{proof}
\begin{thm}
\label{thrm: whenfreefromB}Let \( \Xn ,X_{n+1},\dots ,X_{p}\in (M,\tau ) \)
be self-adjoint non-commutative random variables. Assume that the family \( \Xn  \)
is free from the von~Neumann algebra generated by \( B \) and \( X_{n+1},\dots ,X_{p} \),
and assume that \( B \) is embeddable into the ultrapower of the hyperfinite
II\( _{1} \) factor. Then
\[
\chi (\Xn ,X_{n+1},\dots ,X_{p}|B)=\chi (\Xn )+\chi (X_{n+1},\dots ,X_{p}|B).\]
In particular,
\[
\chi (X_{1},\dots \, X_{n}|B)=\chi (\Xn ).\]
The same statements hold for \( \chi ^{\omega } \).
\end{thm}
\begin{proof}
It is sufficient to show that
\[
\chi (\Xn )+\chi (X_{n+1},\dots ,X_{p}|B)\geq \chi (\Xn ,X_{n+1},\dots ,X_{p}|B).\]
Fix \( m>0 \) and elements \( \Ym \in B \). Let \( \ym \in \Gamma _{R}(\Ym ;k,l,\epsilon ) \);
such an \( m \)-tuple exists because \( B \) can be embedded into the ultrapower
of the hyperfinite II\( _{1} \) factor. By Voiculescu's result in \cite{dvv:improvedrandom},
we have that the ratio
\[
\liminf _{k\to \infty }\frac{\lambda \Gamma _{R}(\Xn ;k,l,\epsilon )\times \Gamma _{R}(X_{n+1},\dots ,X_{p}|y_{1},\dots ,y_{m}[\Ym ];k,l,\epsilon )}{\lambda \Gamma _{R}(\Xn ,X_{n+1},\dots ,X_{p}|\ym [\Ym ];k,l,\epsilon )}\geq 1.\]
Indeed, given a \( \delta >0 \), there is a \( k_{0} \), for all \noun{\( k>k_{0} \)}
and for each choice of an approximant
\[
(x_{n+1},\dots ,x_{p})\in \Gamma _{R}(X_{n+1},\dots ,X_{p}|y_{1},\dots ,y_{m}[\Ym ];k,l,\epsilon ),\]
there exists an open subset
\[
\Gamma _{k,x_{n+1},\dots ,x_{p}}\subset \Gamma _{R}(\Xn ;k,l,\epsilon ),\]
so that
\[
\frac{\lambda \Gamma _{k,x_{n+1},\dots ,x_{p}}}{\Gamma _{R}(X_{1},\dots ,X_{n};k,l,\epsilon )}>1-\delta \]
satisfying
\[
\Gamma _{k,x_{n+1},\dots ,x_{p}}\times \{(x_{n+1},\dots ,x_{p})\}\subset \Gamma _{R}(\Xn ,X_{n+1},\dots ,X_{p}|y_{1},\dots ,y_{m}[\Ym ];k,l,\frac{\epsilon }{2}).\]
Let \( O=O(x_{n+1},\dots ,x_{p}) \) be an open ball of radius \( \epsilon ' \)
for the operator norm on \( M_{k}^{n-p} \), centered at \( x_{n+1},\dots ,x_{p} \).
Then for sufficiently small \( \epsilon ' \) (depending only on \( k \) and
\( l \)), we have
\[
\Gamma _{k,x_{n+1},\dots ,x_{p}}\times O\subset \Gamma _{R}(\Xn ,X_{n+1},\dots ,X_{p}|y_{1},\dots ,y_{m}[\Ym ];k,l,\epsilon ).\]

Let
\[
\Gamma _{k}=\bigcup _{x_{n+1},\dots ,x_{n}}\Gamma _{k,x_{n+1},\dots \, x_{p}}\times O(x_{n+1},\dots ,x_{n}).\]
Then \( \Gamma _{k}\subset \Gamma _{R}(\Xn ,X_{n+1},\dots ,X_{p}|y_{1},\dots ,y_{m}[\Ym ];k,l,\epsilon ), \)
so that
\begin{eqnarray*}
\lambda \Gamma _{R}(\Xn ,X_{n+1},\dots ,X_{p}|y_{1},\dots ,y_{m}[\Ym ];k,l,\epsilon )\geq \lambda \Gamma _{k} &  & \\
\geq \inf _{x_{n+1,\dots \, x_{p}}}\lambda \Gamma _{k,x_{n+1},\dots ,x_{p}}\times \lambda \Gamma _{R}(X_{n+1},\dots ,X_{p}|y_{1},\dots ,y_{m}[\Ym ];k,l,\epsilon ) &  & \\
\geq (1-\delta )\cdot \lambda \Gamma _{R}(\Xn ;k,l,\epsilon )\times \lambda \Gamma _{R}(X_{n+1},\dots ,X_{p}|y_{1},\dots ,y_{m}[\Ym ];k,l,\epsilon ). &  & 
\end{eqnarray*}
The statement of the theorem follows. The proof for \( \chi ^{\omega } \) is
identical.
\end{proof}
We note that the preceding theorem implies that \( \chi (\Xn |B) \) is not
always \( -\infty  \). For example, if \( S_{1},\dots ,S_{n} \) is a free
semicircular family free from a unital von~Neumann algebra \( B \), which can
be embedded into the ultrapower of the hyperfinite II\( _{1} \) factor, then
\( \chi (S_{1},\dots ,S_{n}|B)=\frac{n}{2}\log 2\pi e>-\infty  \).

\begin{cor}
Let \( X_{1},\dots ,X_{p} \) be random variables. Assume that \( \Xn  \) are
free from the algebra generated by \( X_{n+1},\dots ,X_{p} \) and \( B \).
Then
\begin{equation}
\label{eqn:additivity}
\chi (\Xn ,X_{n+1},\dots ,X_{p}|B)=\chi (\Xn |B)+\chi (X_{n+1},\dots ,X_{p}|B).
\end{equation}
 
\end{cor}
\begin{proof}
If \( B \) is embeddable into the ultrapower of the hyperfinite II\( _{1} \)
factor, we have by Theorem \ref{thrm: whenfreefromB} that
\[
\chi (\Xn ,X_{n+1},\dots ,X_{p}|B)=\chi (\Xn )+\chi (X_{n+1},\dots ,X_{p}|B),\]
and also that
\[
\chi (\Xn |B)=\chi (\Xn ).\]
If \( B \) is not embeddable into the ultrapower of the hyperfinite II\( _{1} \)
factor, then the quantities on both sides of (\ref{eqn:additivity}) are equal
to \( -\infty  \).
\end{proof}
The analogy between \( \chi (\Xn |B) \) and \( \chi ^{*}(\Xn :B) \) makes
it tempting to conjecture that (\ref{eqn:additivity}) holds under the weaker
assumption that \( \Xn  \) and \( X_{n+1},\dots ,X_{p} \) are free with amalgamation
over \( B \); however, we were unable to prove this.

\section{Separate change of variables formulas.}

\begin{thm}
\label{thrm: onechangeofvar}Let \( f_{i}:\mathbb {R}\to \mathbb {R} \) be
diffeomorphisms, and let \( \mu _{i} \) be the distribution of \( X_{i} \).
Then
\[
\chi (f_{1}(X_{1}),\dots ,f_{n}(X_{n})|B)=\chi (\Xn |B)+\sum ^{n}_{i=1}\int \int \frac{\log |f_{i}(s)-f_{i}(t)|}{\log |s-t|}d\mu _{i}(s)d\mu _{i}(t),\]
and the same formula holds for \( \chi ^{\omega } \) in place of \( \chi  \).
\end{thm}
\begin{proof}
It is sufficient to prove the statement assuming further that \( f_{i} \) are
identity diffeomorphisms for \( i>1 \); we write \( f=f_{1} \). It is moreover
sufficient to show that given \( \Ym \in B \)
\begin{eqnarray*}
\chi (f_{1}(X_{1}),\dots ,f_{n}(X_{n})|\Ym )\geq \chi (\Xn |\Ym )+ &  & \\
\int \int \frac{\log |f(s)-f(t)|}{\log |s-t|}d\mu _{1}(s)d\mu _{1}(t), &  & 
\end{eqnarray*}
since the reverse inequality follows by replacing \( f \) with its inverse.
It is shown in \cite[Proposition 3.1]{dvv:entropy4} that given \( \delta >0 \),
\( \epsilon >0 \), \( l>0 \), \( R>0 \), there exist \( k_{0}>0 \), \( \epsilon >\epsilon _{0}>0 \),
\( l_{0}>l>0 \), such that for all \( k>k_{0} \), \( 0<\epsilon '<\epsilon _{0} \)
and \( l'>l_{0} \), the determinant of the map
\[
F:(\xn ,\ym )\mapsto (f(x_{1}),x_{2},\dots ,x_{n},\ym )\]
is bounded below by
\[
\exp \left( k^{2}\left[ \int \int \log \frac{|f(s)-f(t)|}{|s-t|}d\mu _{1}(s)d\mu _{1}(t)-\delta \right] \right) \]
for
\[
(\xn ,\ym )\in \Gamma _{R'}(\Xn ,\Ym ;k,l',\epsilon ').\]
Moreover, the image of
\[
\Gamma _{R'}(\Xn ,\Ym ;k,l',\epsilon ')\]
under this map is contained in
\[
\Gamma _{R}(f(X_{1}),\dots ,X_{n},\Ym ;k,l,\epsilon ).\]
Choose \( \ym  \) such that
\begin{eqnarray*}
\log \sup _{z_{1},\dots ,z_{m}}\lambda \Gamma _{R'}(\Xn |z_{1},\dots ,z_{m}[\Ym ];k,l',\epsilon ')- &  & \\
\log \lambda \Gamma _{R'}(\Xn |\ym [\Ym ];k,l;,\epsilon ')<\delta . &  & 
\end{eqnarray*}
Then we have that
\begin{eqnarray*}
\sup _{z_{1},\dots ,z_{m}}\log \lambda \Gamma _{R'}(f(X_{1}),X_{2},\dots ,X_{n}|z_{1},\dots ,z_{m}[\Ym ];k,l,\epsilon )\geq  &  & \\
\log \lambda \Gamma _{R}(f(X_{1}),\dots X_{n}|\ym [\Ym ];k,l,\epsilon )\geq  &  & \\
\log \lambda F\left( \Gamma _{R'}(\Xn |\ym [\Ym ];k,l',\epsilon ')\right) \geq  &  & \\
\log \lambda \Gamma _{R'}(\Xn |\ym [\Ym ];k,l',\epsilon ') &  & \\
+\int \int \frac{\log |f(s)-f(t)|}{\log |s-t|}d\mu _{1}(s)d\mu _{1}(t)-\delta k^{2}\geq  &  & \\
\log \sup _{z_{1},\dots ,z_{m}}\lambda \Gamma _{R'}(\Xn |z_{1},\dots ,z_{m}[\Ym ];k,l',\epsilon ') &  & \\
+\int \int \frac{\log |f(s)-f(t)|}{\log |s-t|}d\mu _{1}(s)d\mu _{1}(t)-2\delta k^{2}. &  & 
\end{eqnarray*}
Taking \( \limsup _{k\to \infty } \) gives us that
\begin{eqnarray*}
\chi _{R}(f(X_{1}),\dots ,X_{n}|\ym [\Ym ];l,\epsilon )\geq  &  & \\
\chi _{R'}(\Xn |\ym [\Ym ];l',\epsilon ') &  & \\
+\int \int \frac{\log |f(s)-f(t)|}{\log |s-t|}d\mu _{1}(s)d\mu _{1}(t), &  & 
\end{eqnarray*}
which implies the theorem. The proof for \( \chi ^{\omega } \) is exactly the
same.
\end{proof}
\begin{prop}
If \( \Vert X_{i}\Vert <R \) for all \( i \), then
\[
\chi _{R}(\Xn |B)=\chi (\Xn |B),\]
and similarly for \( \chi ^{\omega } \).
\end{prop}
The proof is along the lines of that of Theorem~\ref{thrm: onechangeofvar},
using the ideas of a similar Proposition in \cite{dvv:entropy2}, and is therefore
omitted.

\section{General change of variables formula. }

\newcommand{\Fn}{F_{1},\dots ,F_{n}}

\newcommand{\In}{i_{1},\dots ,i_{k}}
Let \( B\subset M \) be a unital subalgebra, and let \( \Fn  \) be non-commutative
power series with coefficients from \( B \); i.e.,
\[
F_{i}(t_{1},\dots ,t_{n})=\sum _{k}\sum _{i_{1},\dots ,i_{k}}b_{\In }^{i,0}t_{i_{1}}b^{i,1}_{\In }\dots t_{i_{k}}b_{\In }^{i,k}.\]
Denote by \( B[t_{1},\dots ,t_{n}] \) the set of all such power series, which
have the property that if \( F\in B[t_{1},\dots ,t_{n}] \) and \( \Xn  \)
are self-adjoint, then \( F(\Xn ) \) is also self-adjoint. Given \( F_{i}\in B[t_{1},\dots ,t_{n}] \)
as above, denote by \( \hat{F}_{i} \) the power series
\[
\hat{F}_{i}(z_{1},\dots ,z_{n})=\sum _{k}\sum _{\In }\prod ^{k}_{j=1}\Vert b_{\In }^{i,j}\Vert y_{i_{1}}\dots y_{i_{k}}.\]
We say that \( (R_{1},\dots ,R_{n}) \) is a mutiradius of convergence of \( F_{i} \),
if it is the multiradius of convergence of \( \hat{F}_{i} \) (as an ordinary
commutative power series). 

Let \( F\in B[t_{1},\dots ,t_{n}] \) be such a power series. Then by the derivative
of \( F \) with respect to \( t_{i} \) we mean the formal power series \( D_{i}F\in B[t_{1},\dots ,t_{n})\otimes B[t_{1},\dots ,t_{n}] \).
Here \( D_{i} \) is defined by the following properties; here we think of \( B[t_{1},\dots ,t_{n}] \)
and \( B[t_{1},\dots ,t_{n})\otimes B[t_{1},\dots ,t_{n}] \) are viewed as
bimodules over the algebra generated by \( B \) and \( t_{1},\dots ,t_{n} \)
using its obvious left and right actions. 

\begin{enumerate}
\item \( D_{i} \) is bilinear over the algebra generated by \( B \) and \( t_{1},\dots ,t_{i-1},t_{i+1},\dots ,t_{n} \);
\item \( D_{i}(t_{i})=1\otimes 1 \); 
\item \( D \) satisfies the Leibniz rule: \( D_{i}(FG)=(D_{i}F)G+F(D_{i}G) \).
\end{enumerate}
As an example,
\[
D_{1}(b_{0}t_{1}b_{1}t_{2}b_{2}t_{1}b_{3}t_{4}b_{4})=b_{0}\otimes b_{1}t_{2}b_{2}t_{1}b_{3}t_{4}b_{4}+b_{0}t_{1}b_{1}t_{2}b_{2}\otimes b_{3}t_{4}b_{4}.\]

Given a family of non-commutative power series \( F_{1},\dots ,F_{n} \) with
a common multiradius of convergence \( (R_{1},\dots ,R_{n}) \), we define for
\( \Xn \in M \), \( \Vert X_{i}\Vert <R_{i} \), its Jacobian at \( \Xn  \),
to be the matrix \( D_{B}F(X_{1},\dots ,X_{n})\in M_{n}\otimes M\otimes M \)
whose \( i,j \)-th entry is equal to \( D_{i}(F_{j})(\Xn ) \). 

Note that if \( B\subset M_{k} \) and \( \xn \in M_{k} \), then \( D_{B}(\xn ) \)
is precisely the Jacobian of the map
\[
(z_{1},\dots ,z_{n})\mapsto \left( F_{1}(z_{1},\dots ,z_{n}),\dots ,F_{n}(z_{1},\dots ,z_{n})\right) ,\]
evaluated at \( \xn  \). 

The proof of the following Theorem is almost identical to the proof of the change
of variables formula given in \cite{dvv:entropy2}, together with the line of
the proof of Theorem~\ref{thrm: onechangeofvar}, and is therefore omitted.

\begin{thm}
\label{thrm: genchangeofvar}Let \( F_{i}\in B[t_{1},\dots ,t_{n}] \), \( i=1,\dots ,n \)
be non-commutative power series with common multiradius of convergence \( (R_{1},\dots ,R_{n}) \).
Assume that there are non-commutative power series \( G_{i} \), \( i=1,\dots ,n \)
with common multiradius of convergence \( (r_{1},\dots ,r_{n}) \), such that
for all \( i=1,\dots ,n \),
\begin{eqnarray*}
F_{i}(G_{1}(t_{1},\dots ,t_{n}),\dots ,G_{n}(t_{1},\dots ,t_{n}))=t_{i}, &  & \\
G_{i}(F_{1}(t_{1},\dots ,t_{n}),\dots ,F_{n}(t_{1},\dots ,t_{n}))=t_{i},. &  & 
\end{eqnarray*}
Assume that for each \( i \), \( \Vert X_{i}\Vert <\min (r_{i},R_{i}) \).
Then
\begin{eqnarray*}
\chi (F_{1}(X_{1},\dots ,X_{n}),\dots ,F_{n}(X_{1},\dots ,X_{n})|B)= &  & \\
\chi (X_{1},\dots ,X_{n}|B)+\Tr \otimes \tau \otimes \tau \left( \log |D_{B}F(\Xn )|\right) . &  & 
\end{eqnarray*}
The same formulas hold for \( \chi ^{\omega } \) in place of \( \chi  \).
\end{thm}
We deduce that ``free Brownian motion'' has a regularizing effect on free
entropy (compare \cite{dvv:entropy5}). The following proposition follows also
from the results of \cite{dvv:entropy3}, but we could not find its exact statement
there.

\begin{prop}
Let \( S_{1},\dots ,S_{n} \) be a free semicircular family, free from the algebra
\( B=W^{*}(\Xn ) \). Assume that \( B \) is embeddable into the ultrapower
of the hyperfinite II\( _{1} \) factor. Then for all \( t>0 \), we have
\[
\chi (X_{1}+\sqrt{t}S_{1},\dots ,X_{n}+\sqrt{t}S_{n})\geq n\log 2\pi et>-\infty .\]
The same estimate holds for \( \chi ^{\omega } \).
\end{prop}
\begin{proof}
By the change of variables formula and Theorem~\ref{thrm: whenfreefromB}, we
have
\begin{eqnarray*}
\chi (X_{1}+\sqrt{t}S_{1},\dots ,X_{n}+\sqrt{t}S_{n})\geq \chi (X_{1}+\sqrt{t}S_{1},\dots ,X_{n}+\sqrt{t}S_{n}|B)= &  & \\
\chi (\sqrt{t}S_{1},\dots ,\sqrt{t}S_{n})=\frac{n}{2}\log 2\pi et. & 
\end{eqnarray*}
 
\end{proof}
\begin{thm}
Let \( P_{1},\dots ,P_{n} \) be non-commutative polynomials in \( n \) variables
with coefficients from \( B \). Assume that \( \chi (\Xn |B)>-\infty  \).
Then
\begin{eqnarray*}
\frac{d}{d\epsilon }\chi (X_{1}+\epsilon P_{1}(X_{1},\dots ,X_{n}),\dots ,X_{n}+\epsilon P_{n}(\Xn )|B)= &  & \\
\sum _{i}\langle J(X_{i}:B\vee W^{*}(X_{1},\dots ,X_{i-1},X_{i+1},\dots ,X_{n}),P_{i}\rangle , &  & 
\end{eqnarray*}
where \( J(X_{i}:B\vee W^{*}(X_{1},\dots ,X_{i-1},X_{i+1},\dots ,X_{n})) \)
is the first-order conjugate variable to \( X_{i} \) with respect to \( B\vee W^{*}(X_{1},\dots ,X_{i-1},X_{i+1},\dots ,X_{n}) \)
(cf. \cite{dvv:entropy5}). The same equality holds for \( \chi ^{\omega } \).
\end{thm}
\begin{proof}
For \( \epsilon >0 \) sufficiently small, the transformation \( F^{\epsilon } \)
defined by
\[
F_{i}^{\epsilon }(X_{1},\dots ,X_{n})=X_{i}+\epsilon P_{i}(X_{1},\dots ,X_{n})\]
is a non-commutative power series in \( \Xn  \) with coefficients from \( B \)
and satisfies the hypothesis of Theorem~\ref{thrm: genchangeofvar}. It follows
that
\begin{eqnarray*}
\chi (X_{1}+\epsilon P_{1}(\Xn ),\dots ,X_{n}+P_{n}(\Xn )|B)=\chi (\Xn |B) &  & \\
+(\Tr \otimes \tau \otimes \tau )\log |D_{B}F^{\epsilon }(\Xn )|. &  & 
\end{eqnarray*}
 Hence the derivative in \( \epsilon  \) of \( \chi (X_{1}+\epsilon P_{1}(\Xn ),\dots ,X_{n}+P_{n}(\Xn )|B) \)
is equal to the derivative of \( (\tau \otimes \tau )\log |D_{B}F^{\epsilon }(\Xn )| \).
Notice that
\[
D_{i}F^{\epsilon }_{j}(\Xn )=\delta _{ij}+\epsilon D_{j}(P_{i})(\Xn ),\]
 so that
\[
D_{B}(F^{\epsilon })(\Xn )=I+\epsilon M,\]
 where \( I \) is the identity matrix and
\[
M_{ij}=D_{i}P_{j}(\Xn ).\]
 Hence
\[
(D_{B}(F^{\epsilon })(\Xn )^{*}D_{B}(F^{\epsilon })(\Xn ))=I+\epsilon (M+M^{*}+\epsilon M^{*}M).\]
 Since \( \log (1+t) \) has a power series expansion around zero, we have that
\[
\frac{1}{2}\log (D_{B}(F^{\epsilon })(\Xn )^{*}D_{B}(F^{\epsilon })(\Xn ))=I+\frac{\epsilon }{2}(M+M^{*})+O(\epsilon ^{2}).\]
 It follows that the desired derivative is equal to
\begin{eqnarray*}
(\tau \otimes \tau )\left( \frac{1}{2}[M+M^{*}]\right) =\frac{1}{2}\sum _{i}\tau \otimes \tau ((D_{i}P_{i})(\Xn )+(D_{i}P_{i})(\Xn )^{*})= &  & \\
\sum _{i}\tau \otimes \tau (D_{i}P_{i})(\Xn ), &  & 
\end{eqnarray*}
since \( F_{i} \) maps self-adjoint variables to self-adjoint variables. Hence
we have, by the definition of the conjugate variable, that
\begin{eqnarray*}
\frac{d}{d\epsilon }\chi (X_{1}+\epsilon P_{1}(\Xn ),\dots ,X_{n}+P_{n}(\Xn )|B)= &  & \\
\sum _{i}\langle J(X_{i}:B\vee W^{*}(X_{1},\dots ,X_{i-1},X_{i+1},\dots ,X_{n}),P_{i}\rangle , & 
\end{eqnarray*}
as claimed.
\end{proof}
Recall (see \cite{dvv:entropy4}) that a function \( \phi (\Xn ) \) is said
to attain a local algebraic maximum at \( \Xn  \) on the set \( S=\{\Xn :\sum _{i}\tau (X_{i}^{2})=n\} \),
if for all non-commutative polynomials \( P_{i}, \) with coefficients from
\( B \) there exist \( \epsilon _{0}>0 \), such that for all \( 0<\epsilon <\epsilon _{0} \),
\[
\phi \left( \frac{X_{1}+\epsilon P_{1}(\Xn )}{\Vert X_{1}+\epsilon P_{1}(\Xn )\Vert _{2}},\dots ,\frac{X_{n}+\epsilon P_{n}(\Xn )}{\Vert X_{n}+\epsilon P_{n}(\Xn )\Vert _{2}}\right) \leq \phi (\Xn ).\]
 Clearly this is a much weaker requirement than saying that \( \phi  \) attains
a maximum on \( S \) at \( \Xn  \).

\begin{prop}
\label{prop: maximality}Let \textbf{\( B \)} be a von~Neumann algebra, embeddable
into an ultrapower of the hyperfinite II\( _{1} \) factor. Then the function
\( (\Xn )\mapsto \chi (\Xn |B) \) attains a local algebraic maximum on the
set \( \{\Xn :\sum _{i}\tau (X_{i}^{2})=n\} \) exactly when \( \Xn  \) are
\( n \) free \( (0,1) \) semicircular variables, free from \( B. \) The same
statement holds for \( \chi ^{\omega } \).
\end{prop}
\begin{proof}
Note that by Corollary~\ref{corr: max}, we have that a global maximum (and
hence a local algebraic maximum) is attained by such a semicircular family.
Assume that the maximum is attained by some family \( \Xn  \). Then \( \chi (\Xn |B) \)
attains a local algebraic maximum at \( \Xn  \). Therefore, we have that for
all non-commutative polynomials \( P_{i} \) with coefficients from \( B \)
\[
\frac{d}{d\epsilon }\chi \left( \left. \frac{X_{1}+\epsilon P_{1}(\Xn )}{\Vert X_{1}+\epsilon P_{1}(\Xn )\Vert _{2}},\dots ,\frac{X_{n}+\epsilon P_{n}(\Xn )}{\Vert X_{n}+\epsilon P_{n}(\Xn )\Vert _{2}}\right| B\right) =0.\]
 But this is equal to
\begin{eqnarray*}
\frac{d}{d\epsilon }\chi (\Xn |B)-\frac{d}{d\epsilon }\sum _{i}\log \Vert X_{i}+\epsilon P_{i}(\Xn )\Vert _{2}= &  & \\
\sum _{i}\langle J(X_{i}:B\vee W^{*}(X_{1},\dots ,X_{i-1},X_{i+1},X_{n})),P_{i}\rangle  &  & \\
-\sum _{i}\langle X_{i},P_{i}(\Xn )\rangle . &  & 
\end{eqnarray*}
It follows that for all non-commutative polynomials \( P_{i} \) with coefficients
from \( B \),
\[
\langle J(X_{i}:B\vee W^{*}(X_{1},\dots ,X_{i-1},X_{i+1},X_{n})),P_{i}\rangle =\langle X_{i},P_{i}(\Xn )\rangle ,\]
which implies that
\[
X_{i}=J(X_{i}:B\vee W^{*}(X_{1},\dots ,X_{i-1},X_{i+1},X_{n}))\]
for all \( i=1,\dots ,n \). But by \cite{dvv:entropy5} and \cite{shlyakht:cpentropy},
this implies that \( \Xn  \) are a free semicircular family, free from \( B \).
\end{proof}
\begin{thm}
Assume that \( B \) is embeddable into an ultrapower of the hyperfinite II\( _{1} \)
factor, and \( \sum _{i=1}^{n}\tau (X_{i}^{2})=n \). Then \( \chi (\Xn |B) \)
attains its maximal value of \( n\log 2\pi e \) if and only if \( \Xn  \)
are a free semicircular family, which is free from \( B \). The same statement
holds for \( \chi ^{\omega } \).
\end{thm}
\begin{proof}
The condition that \( \chi (\Xn |B) \) achieves a local algebraic maximum is
weaker than the condition that it achieves its maximum, so if the maximum is
achieved, the local algebraic maximum is achieved, and Proposition~\ref{prop: maximality}
applies. Conversely, it was shown in Corollary~\ref{corr: max} that the given
number is indeed a maximum.
\end{proof}
We end with the following theorem, whose proof is identical to that of \cite[Proposition 4.3]{dvv:entropy4}.

\begin{thm}
If \( \chi (X|B)=\chi (X)\neq -\infty  \), then \( X \) is free from \( B \).
The same statement holds for \( \chi ^{\omega } \).
\end{thm}

\section{\protect\( \chi (\Xn |B\otimes M_{N})\protect \). }

Let \( X_{ij}^{k} \), \( i,j=1,\dots ,N \), \( k=1,\dots ,n \) be non-commutative
random variables, such that \( X_{ij}=X_{ji}^{*} \), and let \( B \) be a
unital subalgebra of \( (M,\tau ) \) . Then the joint \( * \)-distribution
of the family \( \{X_{ij}^{k}\}\cup B \) completely determines and is completely
determined by, the joint distribution of the matrices
\[
Z_{k}=\left( \begin{array}{cccc}
X_{11}^{k} & X_{12}^{k} & \dots  & X_{1N}^{k}\\
X^{k}_{21} & X_{22}^{k} & \dots  & X_{2N}^{k}\\
\vdots  & \vdots  & \ddots  & \vdots \\
X_{N1}^{k} & X_{N2}^{k} & \dots  & X_{NN}^{k}
\end{array}\right) ,\]
the matrix units \( E_{ij}\in M_{N} \) (matrices whose only non-zero entry
is in the position \( i,j \)) and the algebra \( B\otimes M_{N} \), identified
with those matrices that have entries from \( B \). Therefore, it is natural
to expect a relationship between the free entropy of the entries of the matrix
relative to \( B \) and the free entropy of the matrix relative to the algebra
\( B\otimes M_{N} \) of \( B \)-valued matrices. Such a property is enjoyed
by \( \chi ^{*} \) (introduced by Voiculescu in \cite{dvv:entropy5}; this
property for \( \chi ^{*} \) was proved in \cite{nss:entropy}).

Let us define
\[
Y_{ij}^{k}=\left\{ \begin{array}{cc}
\frac{1}{2}(X_{ij}^{k}+[X_{ij}^{k}]^{*}) & \textrm{if }i\geq j\\
\frac{1}{2\sqrt{-1}}(X_{ij}^{k}-[X_{ij}^{k}]^{*}) & \textrm{if }i<j.
\end{array}\right. \]
For \( \omega  \) a free ultrafilter (i.e., a homomorphism \( \omega :C(\mathbb {N})\to \mathbb {C}) \),
from the algebra of continuous bounded functions on \( \mathbb N \)), and \( n\in \mathbb N \),
define \( n\omega  \) to be the free ultrafilter, which as a homomorphism from
\( C(\mathbb {N}) \) is given by the composition of \( \omega  \) and the
map \( n\cdot f \), given by \( (n\cdot f)(m)=f(nm) \). 

\begin{thm}
Let \( Y_{\{ij\}} \)Let \( \omega  \) be a free ultrafilter. Then
\[
\chi ^{N\omega }(\{Y_{ij}^{k}\}_{i,j,k}|B)=N^{2}\chi ^{\omega }(Z_{1},\dots ,Z_{n}|B\otimes M_{N})-N^{2}\frac{n}{2}\log N.\]
Moreover,
\[
\chi (\{Y_{ij}^{k}\}_{i,j,k})\leq N^{2}\chi (Z_{1},\dots ,Z_{n}|B\otimes M_{N})-N^{2}\frac{n}{2}\log N.\]

\end{thm}
\begin{proof}
Let \( E_{ij} \) be as before. Then
\[
\chi ^{\omega }(Z_{1},\dots ,Z_{n}|B\otimes M_{N})=\inf _{P_{1},\dots ,P_{q}\in B}\chi ^{\omega }(Z_{1},\dots ,Z_{n}|P_{1},\dots ,P_{q},\{E_{ij}\}_{ij}).\]
Here \( E_{ij} \) are not self-adjoint; what we mean by the quantity on the
right is the obvious extension of our quantity to such a non-selfadjoint case. 

We first claim that
\[
\chi ^{N\omega }(\{Y_{ij}^{k}\}_{i,j,k}|B)\geq N^{2}\chi ^{\omega }(Z_{1},\dots ,Z_{n}|B\otimes M_{N})-N^{2}\frac{n}{2}\log N.\]

Assume that \( k=Nk' \). Fix \( \delta >0 \). Choose \( Q_{1},\dots ,Q_{s}\in B\otimes M_{N} \)
and \( R>0 \) so that
\[
\chi _{R}^{\omega }(Z_{1},\dots ,Z_{n}|Q_{1},\dots ,Q_{q},\{E_{ij}\}_{ij})\geq \chi ^{\omega }(Z_{1},\dots \, Z_{n}|B\otimes M_{N})-\delta .\]
 Let \( (e_{ij})\in \Gamma _{R}(\{E_{ij}\}:Z_{1},\dots \, Z_{n},Q_{1},\dots ,Q_{s};k,l,\epsilon ) \).
Then by a suitable choice of \( l \), \( \epsilon  \), \( k \) and \( R \),
we can guarantee that the exists a projection \( p\leq e_{11} \) of rank \( k' \),
\( [e_{11},p]=0 \). Given \( \delta  \), choose \( q_{1},\dots ,q_{S} \)
and \( e_{ij} \) so that
\begin{eqnarray*}
\log \lambda \Gamma _{R}(Z_{1},\dots ,Z_{n}|\{e_{ij}\},\{q_{i}\}[\{E_{ij}\},\{Q_{i}\}];k,l,\epsilon )\geq  &  & \\
\sup _{(f_{ij})\in \Gamma _{R}(\{E_{ij}\}:k,l,\epsilon )}\Gamma _{R}(Z_{1},\dots ,Z_{n}|\{f_{ij}\}[\{E_{ij}\}];k,l,\epsilon )-\delta k^{2} &  & 
\end{eqnarray*}
 Let \( p \) be as before, and identify \( pM_{k}p \) with \( M_{k'} \).
For \( (z_{1},\dots ,z_{n})\in M^{n}_{k} \), let \( z_{ij}^{k}=pe_{1i}z_{k}e_{j1}p \)
and \( q^{r}_{ij}=pe_{1i}q_{r}e_{j1}p \). Denote by \( T=T_{\{e_{ij}\}} \)the
map from \( M_{k}^{n} \) to \( M_{k'}^{nN^{2}} \) given by \( T(z_{1},\dots ,z_{n})=(y^{r}_{ij})_{ijr} \),
where
\[
y^{r}_{ij}=\left\{ \begin{array}{cc}
\frac{1}{2}(z_{ij}^{r}+[z_{ij}^{r}]^{*}) & \textrm{if }i\geq j,\\
\frac{1}{2\sqrt{-1}}(z_{ij}^{r}-[z_{ij}^{r}]^{*}) & \textrm{if }i<j.
\end{array}\right. \]
It follows that for \( l \), \( k \), \( R \) sufficiently large and \( \epsilon  \)
sufficiently small, we can assume that the logarithm of the Jacobian of \( T \)
is at least \( -\delta k^{2} \). Moreover, given \( l', \) \( R' \) and \( \epsilon ' \)
there exist \( l>l' \), \( R>R' \) and \( 0<\epsilon <\epsilon ' \), such
that
\[
T(\Gamma _{R}(Z_{1},\dots ,Z_{n}|\{e_{ij}\},\{q_{i}\}[\{E_{ij}\},\{Q_{i}\}];k,l,\epsilon )\subset \Gamma _{R'}(\{Y_{is}^{r}\}|\{q_{ij}^{r}\}_{ijr}[\{Q_{ij}^{r}\}_{ijr}];k',l',\epsilon ').\]
 It follows that
\begin{eqnarray*}
\log \lambda \Gamma _{R'}(\{Y_{is}^{r}\}|\{q_{ij}^{r}\}_{ijr}[\{Q_{ij}^{r}\}_{ijr}];k',l',\epsilon ')\geq  &  & \\
\log \lambda T(\Gamma _{R}(Z_{1},\dots ,Z_{n}|\{e_{ij}\},\{q_{i}\}[\{E_{ij}\},\{Q_{i}\}];k,l,\epsilon )\geq  &  & \\
\log \lambda \Gamma _{R}(Z_{1},\dots ,Z_{n}|\{e_{ij}\},\{q_{i}\}[\{E_{ij}\},\{Q_{i}\}];k,l,\epsilon )-\delta k^{2}\geq  &  & \\
\sup _{((f_{ij})_{ij},(p_{r})\in \Gamma _{R}(\{E_{ij}\},\{P_{i}\}:k,l,\epsilon )}\log \lambda \Gamma _{R}(Z_{1},\dots ,Z_{n}|\{f_{ij}\},\{p_{i}\}[\{E_{ij}\},\{Q_{i}\}];k,l,\epsilon )-2\delta k^{2}. & 
\end{eqnarray*}
Therefore, remembering that \( k=Nk' \) and taking limits as \( k\to \omega  \),
we get:
\begin{eqnarray*}
\chi ^{N\omega }_{R'}(\{Y_{is}^{r}\}|\{q_{ij}^{r}\}_{ijr}[\{Q_{ij}^{r}\}_{ijr}];l',\epsilon ')= &  & \\
\lim _{k\to \omega }\frac{1}{(k')^{2}}\log \lambda \Gamma _{R'}(\{Y_{is}^{r}\}\{q_{ij}^{r}\}_{ijr}[\{Q_{ij}^{r}\}_{ijr}];k',l',\epsilon ')+\frac{N^{2}n}{2}\log k'= &  & \\
\lim _{k\to \omega }N^{2}\left( \frac{1}{k^{2}}\log \lambda \Gamma _{R'}(\{Y_{is}^{r}\}\{q_{ij}^{r}\}_{ijr}[\{Q_{ij}^{r}\}_{ijr}];k',l',\epsilon ')+\frac{n}{2}\log k-\frac{n}{2}\log N\right) \geq  &  & \\
N^{2}\lim _{k\to \omega }\frac{1}{k^{2}}\sup _{(f_{ij})\in \Gamma _{R}(\{E_{ij}\}:k,l,\epsilon )}\log \lambda \Gamma _{R}(Z_{1},\dots ,Z_{n}|\{f_{ij}\},\{q_{i}\}[\{E_{ij}\},\{Q_{i}\}];k,l,\epsilon ) &  & \\
+N^{2}\frac{n}{2}\log k-N^{2}\frac{n}{2}\log N-2\delta = &  & \\
N^{2}\chi _{R}^{\omega }(Z_{1},\dots ,Z_{n}|\{E_{ij}\},\{Q_{i}\})-N^{2}\frac{n}{2}\log N-2\delta \geq  &  & \\
N^{2}\chi _{R}^{\omega }(Z_{1},\dots ,Z_{n}|B\otimes M_{N};l,\epsilon )-N^{2}\frac{n}{2}\log N-3\delta . &  & 
\end{eqnarray*}
 This implies the claimed inequality. 

Next, we claim that
\[
\chi ^{N\omega }(\{Y_{ij}^{k}\}_{i,j,k}|B)\geq N^{2}\chi ^{\omega }(Z_{1},\dots ,Z_{n}|B\otimes M_{N})-N^{2}\frac{n}{2}\log N.\]

For this, choose \( Q_{1},\dots ,Q_{n}\in B\otimes M_{N} \) and \( R>0 \)
in such a way that
\[
\chi _{R}^{\omega }(Z_{1},\dots ,Z_{n}|Q_{1},\dots \, Q_{n},\{E_{ij}\})\geq \chi ^{\omega }(Z_{1},\dots ,Z_{n}|B\otimes M_{N})-\delta .\]
Set \( Q_{ij}^{r}=E_{1i}Q_{r}E_{j1} \). Choose \( \{q_{ij}^{r}\}\in \Gamma _{R}(\{Q_{ij}^{r}\};m,k,\epsilon ) \).
Assume that
\[
(y_{ij}^{r})\in \Gamma _{R}((Y_{ij}^{r})|\{q_{ij}^{r}\}[\{Q_{ij}^{r}\}];k,l,\epsilon ).\]
Then let \( x_{ij}^{r}=y_{ij}^{r}+\sqrt{-1}y_{ji}^{r} \). Set \( x_{r}\in M_{k}\otimes M_{N} \)
to be \( x_{r}=\sum _{ij}x_{ij}^{r}\otimes E_{ij} \), and put \( q_{r}=\sum q_{ij}^{r}\otimes E_{ij}\in M_{k}\otimes M_{n} \).
Then given \( R'>0 \), \( l'>0 \) and \( \epsilon '>0 \), there exist \( R>R'>0 \),
\( l>l'>0 \) and \( 0<\epsilon <\epsilon ' \), such that \( (z_{1},\dots ,z_{r})\in \Gamma _{R'}(Z_{1},\dots ,Z_{n}|\{E_{ij}\},\{q_{i}\}[\{E_{ij}\},\{Q_{i}\}];nk,l',\epsilon ') \).
Since the map assigning to \( (y_{ij}^{r}) \) the \( n \)-tuple \( (z_{1},\dots ,z_{r}) \)
is measure-preserving, we get that
\begin{eqnarray*}
\sup _{e_{ij},f_{i}}\log \lambda \Gamma _{R'}(Z_{1},\dots ,Z_{n}|\{e_{ij}\},\{f_{i}\}[\{E_{ij}\},\{Q_{i}\}];nk,l',\epsilon ')\geq  &  & \\
\log \lambda \Gamma _{R'}(Z_{1},\dots ,Z_{n}|\{E_{ij}\},\{q_{i}\}[\{E_{ij}\},\{Q_{i}\}];nk,l',\epsilon ')\geq  &  & \\
\log \lambda \Gamma _{R}((Y_{ij}^{r})|\{q_{ij}^{r}\}[Q_{ij}^{r}];k,l,\epsilon ). &  & 
\end{eqnarray*}
This implies our claim.

The proof of the inequality for \( \chi  \) instead of \( \chi ^{\omega } \)
is along the lines of the proof of the second inequality above, and is left
to the reader.
\end{proof}
\bibliographystyle{amsplain}

\providecommand{\bysame}{\leavevmode\hbox to3em{\hrulefill}\thinspace}

\end{document}